\documentclass[a4paper]{amsart}
\usepackage{amsmath,amssymb,amsfonts}
\usepackage{colordvi}
\usepackage[all]{xy}

\def\mpoint{\;.}
\def\mvirg{\;,}

\def\sou{\underline}

\def\mpn{\medskip\par\noindent}

\def\mmpn{\vskip 1em minus 1em\par\noindent}

\def\sp{\bigskip\par}

\def\smp{\smallskip\par}

\def\CA{{\mathcal A}}
\def\CB{{\mathcal B}}

\def\CD{{\mathcal D}}

\def\CF{{\mathcal F}}

\def\CL{{\mathcal L}}

\def\CZ{{\mathcal Z}}

\def\Im{\operatorname{Im}\nolimits}
\def\Id{\operatorname{id}\nolimits}

\def\Hom{\operatorname{Hom}\nolimits}
\def\End{\operatorname{End}\nolimits}

\def\tot{\operatorname{tot}\nolimits}

\def\op{^{op}}

\def\meet{\wedge}

\def\N{\mathbb{N}}

\newcommand{\edge}[2]{\xymatrix{#1\ar@{->-}[r]&#2}}
\def\marc[#1]{\ar@{-}[#1]|(.4){\object@{<}}}
\def\mard[#1]{\ar@{-}[#1]|(.5){\object@{>}}}
\def\marb[#1]{\ar@{-}[#1]|{\object+{  }}}
\newcommand{\fleche}[2]{\xymatrix@C=4ex{*!U(0.2){#1\;}&*!U(0.5){\;#2}\marc[l]}}
\newcommand{\flecheb}[2]{\xymatrix@C=4ex{*!U(0.2){#1\;}&*!U(0.1){\;#2}\marc[l]}}

\def\pf{\par\bigskip\noindent{\bf Proof~: }}
\def\endpf{~\hfill\rlap{\hspace{-1ex}\raisebox{.5ex}{\framebox[1ex]{}}\sp}\bigskip\pagebreak[3]}

\makeatletter
\renewenvironment{enumerate}{\ifnum \@enumdepth >3 \@toodeep\else
       \advance\@enumdepth \@ne
       \edef\@enumctr{enum\romannumeral\the\@enumdepth}\list
       {\csname  label\@enumctr\endcsname}{\setlength{\topsep}{1ex}
 \setlength{\itemsep}{0 pt}\usecounter
         {\@enumctr}\def\makelabel##1{\hss\llap{##1}}}\fi}{\endlist}

\def\@seccntformat#1{\csname the#1\endcsname.\quad}

\def\section{\pagebreak[3]\setcounter{prop}{0}\setcounter{equation}{0}\@startsection{section}{1}{\z@}{4ex plus  6ex}{2ex}{\center\reset@font \large\bf}}

\makeatother

\def\theprop{\thesection.\arabic{prop}}

\renewenvironment{equation}{\refstepcounter{subsection}\refstepcounter
{prop}$$}{\leqno{\bf (\theprop)}$$}

\newenvironment{enonce}[1]{\pagebreak[3]\refstepcounter{prop}\mmpn
{{\bf  \thesection.\arabic{prop}.\ #1.}}\begin{it} }{\end{it}\smp}
\def\thesection{\arabic{section}}

\newcommand{\result}[1]{\begin{enonce}{#1}}
\newcommand{\fresult}{\end{enonce}}

\begin{document}

\title[A central idempotent]
{A central idempotent in the endomorphism algebra of a finite lattice} 

\author{Serge Bouc}
\author{Jacques Th\'evenaz}
\date\today



\begin{abstract}
We give a direct construction of a specific idempotent in the endomorphism algebra of a finite lattice $T$.
This idempotent is associated with all possible sublattices of $T$ which are total orders.
\end{abstract}

\maketitle


\section{Introduction}

\noindent
Let $T$ be a finite lattice and let $k$ be a commutative ring.
The set of all $k$-linear combinations of join-morphisms from $T$ to~$T$ is a $k$-algebra $\End_{k\CL}(T)$ which plays an important role in our work on correspondence functors \cite{BT1, BT2, BT3}.
This algebra is also of independent interest from a purely combinatorial point of view because it reflects the structure of~$T$ in an algebraic fashion.
Here $\CL$ refers to the category of finite lattices, defined in Section~\ref{Section-lattices}, and $k\CL$ is its $k$-linearization.\par

We introduced in~\cite{BT2} an idempotent $e_T^{\tot} \in \End_{k\CL}(T)$ which is associated with all possible subsets of~$T$ which are totally ordered.
We proved that $e_T^{\tot}$ is central and that $e_T^{\tot}\End_{k\CL}(T)$ is isomorphic to a product of matrix algebras.
Unfortunately, the definition of $e_T^{\tot}$ relies on some rather cumbersome constructions.\par

In the present paper, we give a new point of view for this idempotent.
We express it by means of a much easier and explicit formula, which also has the advantage of allowing for computer calculations.
This formula does not depend on our previous work, but of course the proof that the result coincides with the idempotent $e_T^{\tot}$ relies on~\cite{BT2}.


\section{Finite lattices} \label{Section-lattices}

\noindent
In this section, we recall the basic facts we need about the category of finite lattices.
For the rest of this paper, $T$ denotes a finite lattice.
We write $\leq$ for its partial order (or $\leq_T$ when necessary), $\vee$ for its join, $\meet$ for its meet, $\hat0=\hat0_T$ for its least element, and $\hat1=\hat1_T$ for its greatest element.
Recall that an empty join is equal to~$\hat0$, while an empty meet is equal to~$\hat1$.

If $T'$ is another finite lattice, a {\em join-morphism} $\varphi: T \to T'$ is a map such that, for any subset $X\subseteq T$, we have
$$\varphi\big(\bigvee_{x\in X} x\big) = \bigvee_{x\in X} \varphi(x) \mpoint$$
The case $X=\emptyset$ yields the property $\varphi(\hat0)=\hat0$.
Recall that, because $T$ is finite, the meet is uniquely determined by the join thanks to the finite expression
$$x\meet y = \bigvee_{\substack{a\in T \\ a\leq x, a\leq y}} a \mpoint$$
However, a join-morphism need not respect the meet, and in particular need not map $\hat1$ to~$\hat1$.\par

We let $\CL$ be the category whose objects are the finite lattices and morphisms are the join-morphisms.
We let $k\CL$ be the $k$-linearization of~$\CL$.
Its objects are again the finite lattices and $\Hom_{k\CL}(T,T')$ is the free $k$-module with basis $\Hom_\CL(T,T')$.
Composition in~$k\CL$ is the $k$-bilinear extension of composition in~$\CL$.
In particular, $\End_{k\CL}(T)=\Hom_{k\CL}(T,T)$ is a $k$-algebra with respect to composition and its $k$-basis is the monoid $\End_\CL(T)$ of all join-endomorphisms of~$T$.\par

The opposite partial order on a finite lattice~$T$ yields the {\em opposite lattice} $T\op$, swapping the role of $\vee$ and $\meet$,
and with $\hat0_{T\op}=\hat1_T$ and $\hat1_{T\op}=\hat0_T$.
Associated with a join-morphism $\varphi: T \to T'$, there is its opposite
$$\varphi\op : T'{\op} \longrightarrow T\op \,, \qquad \varphi\op(t')= \bigvee_{\substack{t\in T \\ \varphi(t) \leq t'}} t \mpoint$$

\result{Lemma} \label{opposite-morphism}
Let $\varphi: T \to T'$ be a join-morphism between two finite lattices.
\begin{enumerate}
\item $\varphi\op:T'{\op} \to T\op$ is a join-morphism (that is, a meet-morphism $T'\to T$).
\item $(\varphi\op)\op=\varphi$.
\item If $\varphi$ is surjective, then 
$$\varphi\op(t')= \bigvee_{\substack{t\in T \\ \varphi(t) = t'}} t \;=\;\sup\{t\in T \mid \varphi(t) = t' \} \mpoint$$
\end{enumerate}
\fresult

\pf See Lemma~8.1 in~\cite{BT2}.
\endpf

Recall that a {\em chain} in~$T$ is a totally ordered subset of~$T$.
If $n\in\N$, we write $\sou n=\{0,1,\ldots,n\}$, a totally ordered lattice with $\hat0_{\sou n}=0$ and $\hat1_{\sou n}=n$.
It is straightforward to see that a join-morphism $\varphi:\sou n\to T$ is simply an order-preserving map such that $\varphi(\hat0)=\hat0$.
Therefore, an injective join-morphism $\varphi:\sou n\to T$ corresponds to a chain
$A=\{a_0,a_1,a_2,\ldots,a_n\}$ in~$T$ such that $a_0=\hat0$, where $a_i=\varphi(i)$.
We let $\CA_{T,n}$ be the set of all chains of size $(n+1)$ in~$T$ whose least element is $a_0=\hat0$.
If $n=0$, there is just one element in $\CA_{T,0}$, namely the chain consisting of $\hat0=a_0$.
\par

Similarly, a surjective join-morphism $\pi: T\to \sou n$ corresponds to a chain
$B=\{b_0,b_1,\ldots,b_{n-1}, b_n\}$ in~$T$ such that $b_n=\hat1$, where
$$b_i=\pi\op(i)=\sup\{t\in T \mid \pi(t) = i\} \mpoint$$
We let $\CB_{T,n}$ be the set of all chains of size $(n+1)$ in~$T$ whose greatest element is $b_n=\hat1$.\par

The set
$$\CA_T:=\bigcup_{n\geq0} \CA_{T,n}$$
is partially ordered by inclusion.
It has no greatest element (unless $T$ is totally ordered) and we let $\infty$ be an additional element,
larger than any $A\in \CA_T$. This allows us to consider the M\"obius function $\mu(A,\infty)$,
or in other words the reduced Euler characteristic $\widetilde\chi \big(\, ]A,\infty[ \,\big)$ of the interval $]A,\infty[$ of all chains containing~$A$.
For later use, we now show that this M\"obius function can be expressed in terms of the M\"obius function of~$T$.

\result{Lemma} \label{Mobius}
Let $A=\{a_0,a_1,a_2,\ldots,a_n\}$ be an element of~$\CA_{T,n}$. Then
$$\mu(A,\infty) = (-1)^{n+1}\prod_{k=1}^n \mu(a_{k-1},a_k) \mvirg$$
where $\mu(a_{k-1},a_k)$ denotes the M\"obius function of the interval $]a_{k-1},a_k[$ in~$T$.
\fresult

\pf
For any poset $X$, let $s_i(X)$ be the number of chains of cardinality~$i$ in~$X$.
For $i=0$, there is the empty chain, so $s_0(X)=1$.
It is well-known that
$$\mu(A,\infty)=\widetilde\chi \big(\,  ]A,\infty[ \,\big)=\sum_{i\geq0} (-1)^{i-1} s_i(\, ]A,\infty[ \,) \mpoint$$
The sign is $(-1)^{i-1}$ because a chain of cardinality~$i$ is an $(i{-}1)$-simplex.
Now if $A'$ is a chain with $A\subseteq A'$, then $A'$ is obtained from~$A$ by inserting a chain in each interval $]a_{k-1},a_k[$ independently. Therefore
$$s_i(\, ]A,\infty[ \,) = \displaystyle \sum_{\substack{i_1,\ldots,i_n\geq0 \\ i_1+\ldots+i_n=i}} \prod_{k=1}^n s_{i_k}(\, ]a_{k-1},a_k[ \,)$$
and it follows that
\begin{eqnarray*}
\mu(A,\infty)&=& \displaystyle \sum_{i_1,\ldots,i_n\geq0} (-1)^{i_1+\ldots +i_n-1} \, \prod_{k=1}^n s_{i_k}(\, ]a_{k-1},a_k[ \,) \\
&=& (-1) \prod_{k=1}^n \Big( \sum_{i_k\geq 0} (-1)^{i_k} s_{i_k}(\, ]a_{k-1},a_k[ \Big) \\
&=& (-1)^{n+1} \prod_{k=1}^n \Big( \sum_{i_k\geq 0} (-1)^{i_k-1} s_{i_k}(\, ]a_{k-1},a_k[ \Big) \\
&=& (-1)^{n+1}\prod_{k=1}^n \mu(a_{k-1},a_k) \mvirg
\end{eqnarray*}
as was to be shown.
\endpf

There is one case when the M\"obius function vanishes.

\result{Lemma} \label{contraction}
Let $A=\{a_0,a_1,a_2,\ldots,a_n\}$ be an element of~$\CA_{T,n}$.
If $a_n<\hat1$, then $\mu(A,\infty)=0$.
\fresult

\pf
For any chain $B=\{b_0,b_1,b_2,\ldots,b_m\}$, we let $\overline B= \{b_0,b_1,b_2,\ldots,b_m, \hat1\}$ if $b_m<\hat1$
and $\overline B=B$ if $b_m=\hat1$.
The poset $]A,\infty[$ is conically contractible in the sense of Quillen (see 1.5 in~\cite{Qu}), via the contraction
$$B \leq \overline B \geq \overline A$$
and it follows that $\mu(A,\infty)=0$.
\endpf

Because of this lemma, we shall only be interested in the subset $\CZ_T\subseteq \CA_T$ consisting of all chains $A$ whose greatest element is~$\hat1$ (and least element~$\hat0$), i.e. such that $\overline A=A$. Thus for any chain $A=\{a_0,a_1,a_2,\ldots,a_n\}$ in~$\CZ_T$, we have
$$\hat0=a_0<a_1<\ldots<a_n=\hat1 \mpoint$$


\section{The idempotent corresponding to total orders} \label{Section-idempotent}

\noindent
In this section, we consider a two sided-ideal $\End_{k\CL}^{\tot}(T)$ of the $k$-algebra $\End_{k\CL}(T)$, corresponding to total orders.
This ideal was considered in Section~10 of~\cite{BT2} and it has a central identity element $e_T^{\tot}\in \End_{k\CL}^{\tot}(T)$.
Our main purpose is to prove that $e_T^{\tot}$ can be expressed by a much simpler formula and to prove it by direct combinatorial arguments.\par

We define $\End_\CL^{\tot}(T)$ to be the subset of~$\End_\CL(T)$ consisting of all join-morphisms $\alpha:T\to T$ such that the image $\alpha(T)$ is a totally ordered subset of~$T$.
We let $\End_{k\CL}^{\tot}(T)$ be the $k$-linear span of~$\End_\CL^{\tot}(T)$ in~$\End_{k\CL}(T)$.

\result{Lemma} \label{ideal}
$\End_{k\CL}^{\tot}(T)$ is a two-sided ideal of $\End_{k\CL}(T)$.
\fresult

\pf
Let $\alpha\in \End_\CL^{\tot}(T)$ and $\varphi \in \End_\CL(T)$. It is clear that the image of $\alpha\varphi$ is totally ordered,
so $\alpha\varphi\in \End_\CL^{\tot}(T)$.
On the other hand, the totally ordered subset $\alpha(T)$ is mapped by~$\varphi$ to a totally ordered subset, so $\varphi\alpha\in \End_\CL^{\tot}(T)$.
The result follows by considering $k$-linear combinations.
\endpf

The following result is Theorem~10.8 of~\cite{BT2} and is the starting point of the present work.

\result{Theorem} \label{decomposition}
There is a subalgebra $\CD$ of $\End_{k\CL}(T)$ such that
$$\End_{k\CL}(T) = \End_{k\CL}^{\tot}(T) \times \CD \mvirg$$
(where $\End_{k\CL}^{\tot}(T)$ is identified with $\End_{k\CL}^{\tot}(T)\times\{0\}$ and $\CD$ with $\{0\}\times\CD$, as usual).
\fresult

We let $e_T^{\tot}$ be the identity element of the factor~$\End_{k\CL}^{\tot}(T)$.
This is a central idempotent of~$\End_{k\CL}(T)$. The identity element $\Id_T\in \End_{k\CL}(T)$ decomposes as
$$\Id_T= e_T^{\tot} + (\Id_T- e_T^{\tot}) \mvirg$$
and $\Id_T- e_T^{\tot}\in\CD$.
The formula for $e_T^{\tot}$ given in Theorem~10.8 of~\cite{BT2} comes from rather elaborate constructions, which we revisit in Section~\ref{Section-original-approach} below. We now give an alternative formula for~$e_T^{\tot}$.\par

For any $B\in \CZ_T$, we define
$$\alpha_B:T\to T \,, \qquad \alpha_B(t):=\min\{b\in B \mid b\geq t \} \mpoint$$
Note that the set $\{b\in B \mid b\geq t \}$ is nonempty because $\hat1\in B$ (using our assumption that $B\in\CZ_T$).
The image of~$\alpha_B$ is equal to~$B$, hence totally ordered.
It follows easily that $\alpha_B$ is a join-morphism.
Thus $\alpha_B\in \End_\CL^{\tot}(T)$ and it is moreover clear that $\alpha_B(t)\geq t$ for any $t\in T$ and that $\alpha^2=\alpha$, because $\alpha(b)=b$ for any $b\in B$.

\result{Theorem} \label{main-theorem}
$e_T^{\tot} = \displaystyle - \sum_{B\in\CZ_T} \mu(B,\infty) \,\alpha_B$.
\fresult

\result{Remarks}
{\rm 
(a) We could as well define $\alpha_B$ for $B\in\CA_T$ and sum over all $B\in\CA_T$, but since $\mu(B,\infty)=0$ whenever $B\in \CA_T-\CZ_T$ by Lemma~\ref{contraction}, we see that we only need to consider a sum indexed by~$\CZ_T$.\par

(b) The sum could be restricted further to all $B\in\CZ_T$ such that the lattice $[b_{k-1},b_k]$ is complemented for each~$k=1,\ldots,n$ (where $B=\{b_0,b_1,\ldots,b_n\}$ and $b_0=\hat0$, $b_n=\hat1$), because
$$\mu(B,\infty) = (-1)^{n+1}\prod_{k=1}^n \mu(b_{k-1},b_k)$$
by Lemma~\ref{Mobius} and $\mu(b_{k-1},b_k)=0$ whenever the lattice $[b_{k-1},b_k]$ is not complemented, by Crapo's formula (see Exercice 92 of Chapter 3 in~\cite{St}).\par

(c) The sum could also be indexed by all endomorphisms $\alpha\in\End_\CL^{\tot}(T)$
satisfying $\alpha\geq\Id$ and $\alpha^2=\alpha$,
because the latter two conditions imply that $\alpha=\alpha_B$ where $B$ is the image of~$\alpha$ (which is totally ordered).
}
\fresult

\bigskip\noindent
{\bf Proof of Theorem~\ref{main-theorem}~: }
Let $e:=\displaystyle - \sum_{B\in\CZ_T} \mu(B,\infty) \, \alpha_B$.
We claim that it suffices to prove that
\begin{equation} \label{e-psi}
e \psi=\psi \,,\qquad \forall \,\psi\in \End_\CL^{\tot}(T) \mpoint
\end{equation}
If (\ref{e-psi}) holds, then
$$e \,e_T^{\tot}=e_T^{\tot}$$
because $e_T^{\tot}$ belongs to~$\End_{k\CL}^{\tot}(T)$ and is therefore a $k$-linear combination of morphisms $\psi\in \End_\CL^{\tot}(T)$.
On the other hand, $e\, e_T^{\tot}=e$ because $e$ belongs to~$\End_{k\CL}^{\tot}(T)$ and $e_T^{\tot}$ is its identity element.
Thus $e_T^{\tot}=e\,e_T^{\tot}=e$, as required.\par

In order to establish~(\ref{e-psi}), we prove more generally that $e \psi=\psi$
for any map $\psi: S\to T$ such that $\Im(\psi)$ belongs to~$\CA_T$ where $S$ is some finite set
(i.e. $\Im(\psi)$ is a totally ordered subset of~$T$ starting with~$\hat0$).
Letting $X=\Im(\psi)$, we decompose $\psi$ as the composite of a surjection $S\to X$ followed by the inclusion map $i_X:X\to T$.
It suffices to prove that $e \,i_X=i_X$ for any chain~$X$ in~$T$ starting with~$\hat0$.
Now we have
\begin{equation} \label{e-iX}
e \,i_X= - \sum_{B\in\CZ_T} \mu(B,\infty) \,\alpha_B i_X
= \sum_{\varphi:X\to T}\big( - \sum_{\substack{B\in\CZ_T \\ \alpha_B i_X=\varphi}} \mu(B,\infty) \big) \,\varphi \mpoint
\end{equation}
By the definition of~$\alpha_B$, the equation $\alpha_B i_X=\varphi$ means that, for any $x\in X$,
the element $\varphi(x)$ is the least element of~$B$ such that $\varphi(x)\geq x$.
In other words, the condition $\alpha_B i_X=\varphi$ is equivalent to
\begin{equation} \label{condition1}
[x,\varphi(x)] \cap B = \{\varphi(x)\} \,, \;\forall \, x\in X \mpoint
\end{equation}
Any function $\varphi:X\to T$ appearing in the sum~(\ref{e-iX}) must satisfy the following 3 conditions~:
\begin{enumerate}
\item $\varphi$ is order-preserving and $\varphi(\hat0)=\hat0$ (that is, $\varphi$ is a join-morphism).
\item $\varphi(x)\geq x$, for all $x\in X$.
\item If $x,y\in X$ satisfy $x\leq y\leq \varphi(x)$, then $\varphi(y)= \varphi(x)$.
\end{enumerate}
In order to prove this, we note that the coefficient of~$\varphi$ in~(\ref{e-iX}) is nonzero only if there exists at least one $B\in\CZ_T$ such that $\alpha_B i_X=\varphi$.
Condition (a) follows from the fact that both $i_X$ and~$\alpha_B$ are order-preserving and map $\hat0$ to~$\hat0$, hence $\varphi=\alpha_B i_X$ has the same properties.
Condition~(b) is clear because $\alpha_B\geq\Id$.
For condition~(c), note that the only element of~$[x,\varphi(x)]\cap B$ is $\varphi(x)$,
so the definition of~$\alpha_B$ yields $\alpha_B(y)=\varphi(x)$,
that is, $\varphi(y)= \varphi(x)$.\par

Now we prove that, if $\varphi$ satisfies (b) and~(c), then (\ref{condition1}) is equivalent to
\begin{equation} \label{condition2}
\big( \bigcup_{x\in X} [x,\varphi(x)] \big) \cap B =\varphi(X) \mpoint
\end{equation}
It is clear that (\ref{condition1}) implies (\ref{condition2}).
Assume now (\ref{condition2}) and let $x\in X$.
Notice that $[x,\varphi(x)] \cap B$ is nonempty because $\varphi(x)\in\varphi(X)$, hence $\varphi(x)\in B$ by~(\ref{condition2}),
and so $\varphi(x)\in[x,\varphi(x)] \cap B$.
For any $b\in [x,\varphi(x)] \cap B$, we have $b\in\varphi(X)$ by~(\ref{condition2}), that is, $b=\varphi(y)$ for some $y\in X$.
Since $X$ is totally ordered, we have either $x\leq y\leq \varphi(y)=b\leq \varphi(x)$, hence $\varphi(y)=\varphi(x)$ by~(c),
or $y\leq x\leq b=\varphi(y)$, hence $\varphi(x)=\varphi(y)$ by~(c) again.
Therefore $b=\varphi(x)$, showing that $[x,\varphi(x)] \cap B=\{\varphi(x)\}$. This proves that (\ref{condition2}) implies~(\ref{condition1}).\par

We now fix a map $\varphi:X\to T$ satisfying (a), (b), (c), and we set
$$C=\bigcup_{x\in X} [x,\varphi(x)] \big) \qquad \text{ and } \qquad D=\varphi(X) \mvirg$$
so that (\ref{condition2}) becomes $C\cap B= D$.
Since $\varphi$ is order-preserving and $\varphi(\hat0)=\hat0$, we have $D\in\CA_T$.
Let $\overline C=C\cup\{\hat1\}$ and $\overline D=D\cup\{\hat1\}$, so that $\overline D\in\CZ_T$.
Clearly, $\hat1\notin C$ if and only if $\hat1\notin D$,
and therefore the condition $C\cap B= D$ is equivalent to $\overline C\cap B= \overline D$.
It follows that the coefficient of~$\varphi$ in~(\ref{e-iX}) is equal to
$$-\sum_{\substack{B\in\CZ_T \\ \alpha_B i_X=\varphi}} \mu(B,\infty)
= - \sum_{\substack{B\in\CZ_T \\ C\cap B= D}} \mu(B,\infty)
= - \sum_{\substack{B\in\CZ_T \\ \overline C\cap B= \overline D}} \mu(B,\infty) \mvirg$$
because the condition $\alpha_B i_X=\varphi$ is equivalent to~(\ref{condition2}) by the discussion above.
\par

By the defining property of the M\"obius function, we obtain
$$-\sum_{\substack{B\in\CZ_T \\ \overline C\cap B= \overline D}} \mu(B,\infty)
=\sum_{\substack{B,A\in\CZ_T \\ B\subseteq A \\ \overline C\cap B= \overline D}} \mu(B,A)
=\sum_{\substack{A\in\CZ_T \\ \overline D\subseteq A }}
\sum_{\substack{B\in\CZ_T \\ B\subseteq A \\ \overline C\cap B= \overline D}} \mu(B,A) \mpoint$$
For a fixed $A\in\CZ_T$, the chain $B$ runs over the interval $[\overline D,A]$
with the additional condition $\overline C\cap B= \overline D$,
which can also be written $(\overline C\cap A)\cap B= \overline D$ because $B\subseteq A$.
By a well-known property of the M\"obius function (Corollary~3.9.3 in \cite{St}), the corresponding sum
$$\sum_{\substack{B\in [\overline D,A] \\ (\overline C\cap A)\cap B= \overline D}} \mu(B,A)$$
is zero, provided the fixed element $\overline C\cap A$ is not equal to the top element~$A$.
If otherwise $\overline C\cap A=A$, then $A\subseteq \overline C$ and $B=\overline C\cap B=\overline D$,
so that the sum over~$B$ has the single term $\mu(\overline D,A)$ for $B=\overline D$.\par

Going back to the coefficient of~$\varphi$ in~(\ref{e-iX}), we obtain
$$ - \sum_{\substack{B\in\CZ_T \\ \alpha_B i_X=\varphi}} \mu(B,\infty)
= \sum_{\substack{A\in\CZ_T \\ \overline D\subseteq A \subseteq \overline C}} \mu(\overline D,A)
= \sum_{\substack{A\in\CZ_{\overline C} \\ \overline D\subseteq A}} \mu(\overline D,A)
= - \mu(\overline D,\infty)$$
where the latter symbol $\infty$ denotes a top element added to the poset $\CZ_{\overline C}$ (consisting of all chains in~$\overline C$ having least element~$\hat0$ and greatest element~$\hat1$).\par 

Recall that $\varphi\geq i_X$ by condition~(b).
We now assume that $\varphi> i_X$ and we want to prove that $\mu(\overline D,\infty)=0$.
Let $y\in X$ be minimal such that $\varphi(y)>y$.
We claim that, for any $A\in \CZ_{\overline C}$, the union $A\cup\{y\}$ is totally ordered.
We have to prove that any $a\in A$ is comparable with~$y$.
Since $a\in \overline C$, either $a=\hat1$ and then we are done because $y\leq \hat1$, or there exists $x\in X$ such that $a\in [x,\varphi(x)]$.
If $y\leq x$, then $y\leq a$ and we are done again.
We can assume now that $y\not\leq x$, hence $x<y$ because $X$ is totally ordered.
By minimality of~$y$, we must have $\varphi(x)=x$, hence $[x,\varphi(x)]=\{x\}$ and $a=x$.
It follows that $a<y$. This completes the proof that $A\cup\{y\}$ is totally ordered.\par

We claim now that $y$ does not belong to~$D=\varphi(X)$.
Otherwise $y=\varphi(z)$ for some $z\in X$.
If we had $z=y$, we would obtain $\varphi(y)=y$, contrary to the choice of~$y$.
It follows that the relation $z\leq\varphi(z)=y$ must be a strict inequality $z<\varphi(z)$.
This contradicts the minimality of~$y$ and proves the claim.
Moreover, $y\notin \overline D$ because $y<\varphi(y)$, hence $y\neq\hat1$.
Consequently, the poset $]\overline D,\infty[$ is conically contractible (see 1.5 in~\cite{Qu})
via the contraction
$$A \leq A\cup\{y\} \geq \overline D\cup\{y\}$$
and it follows that $\mu(\overline D,\infty)=0$.\par

This shows that the coefficient of~$\varphi$ in~(\ref{e-iX}) is zero whenever $\varphi> i_X$.
Therefore we are left with a single term for $\varphi=i_X$, namely
$$e \,i_X= -\mu(\overline D, \infty) \, i_X \mpoint$$
But for $\varphi=i_X$, we have
$$C=\bigcup_{x\in X} [x,\varphi(x)]=\bigcup_{x\in X} \{x\} = X = i_X(X) = D \mvirg$$
and consequently the only chain in~$\overline C$ containing~$\overline D$ is $\overline D$ itself.
In other words $]\overline D, \infty[=\emptyset$ and $\mu(\overline D, \infty)=-1$.
The required equality $e \,i_X=  i_X$ follows and this completes the proof of Theorem~\ref{main-theorem}.
\endpf

\result{Remark} \label{remark}
{\rm 
It is easy to prove directly that the expression
$$e=-\sum_{B\in\CZ_T} \mu(B,\infty)\, \alpha_B$$
is idempotent, because (\ref{e-psi}) implies that
$e\alpha_B=\alpha_B$ for any $B\in\CZ_T$, hence 
$$e^2=e\big(-\sum_{B\in\CZ_T} \mu(B,\infty)\, \alpha_B\big) = -\sum_{B\in\CZ_T} \mu(B,\infty)\, \alpha_B =e \mpoint$$
However, the proof that this idempotent is central is more elaborate and appears in Theorem~10.8 of~\cite{BT2}.
}
\fresult


\section{The original approach to the idempotent} \label{Section-original-approach}

\noindent
The idempotent $e_T^{\tot}$ was defined in Section~10 of~\cite{BT2} by an explicit formula.
Using this formula, we want to prove that $e_T^{\tot}$ satisfies the equation of Theorem~\ref{main-theorem}.
In other words, we are going to provide a second proof of that theorem, based on the original approach of~\cite{BT2}.
We first need to define the notation.\par

For any $n\in\N$, we use the set $\CB_{T,n}$ of all chains $B=\{b_0,b_1,\ldots,b_n\}$ in~$T$ whose greatest element is $b_n=\hat1$.
We have seen in Section~\ref{Section-lattices} that the set $\CB_{T,n}$ parametrizes the set of surjective join-morphism $\pi: T\to \sou n$ via the rule
$$b_i=\pi\op(i)=\sup\{t\in T \mid \pi(t) = i\} \mpoint$$
Instead of~$\sou n$, it will be convenient to use a totally ordered lattice~$P$ of cardinality $n+1$, that is, a lattice isomorphic to~$\sou n$,
and to define $r(p)=\sup\{q\in P\mid q<p\}$, for any $p\in P-\{\hat0\}$.
With this notation, a surjective join-morphism $\pi: T\to P$ corresponds to a chain $B=\{b_p\mid p\in P\}$ defined by
$$b_p=\pi\op(p)=\sup\{t\in T \mid \pi(t) = p\} \mvirg$$
and satisfying $b_p<b_q$ whenever $p<q$.
We write $\pi^B:T\to P$ for the surjective join-morphism corresponding to the chain $B\in\CB_{T,n}$.
Then $\pi^B(t)=\hat0$ if $t\leq b_0$ and otherwise we recall the rule
$$\pi^B(t)=p \qquad \text{ if } \;t\leq b_p \;\text{ and } \; t\not\leq b_{r(p)} \mpoint$$

For any given $B\in\CB_{T,n}$, we choose an element $a_p\in [b_{r(p)},b_p]$ for each $p\in P-\{\hat0\}$.
This defines a family $A=(a_p)_{p\in P-\{\hat0\}}$ of elements of~$T$.
We let $\CF_B$ be the set of all families $A=(a_p)_{p\in P-\{\hat0\}}$ of elements of~$T$ such that
$a_p\in [b_{r(p)},b_p]$ for every $p\in P-\{\hat0\}$.
If $A\in\CF_B$, we also set $a_{\hat 0}=\hat 0$ and we define
$$j_A^B: P \longrightarrow T \,, \qquad j_A^B(p)= a_p \mpoint$$
Clearly $j_A^B$ is order-preserving (because if $p<q$ in~$P$,
then $p\leq r(q)$, hence $a_p\leq b_p\leq b_{r(q)}\leq a_q$),
and it also maps $\hat0$ to~$\hat0$. Therefore $j_A^B$ is a join-morphism.

Now let $B^-=\{b_{r(p)} \mid p\in P-\{\hat0\} \}$ and for any $A\in\CF_B$, write
$$\mu(B^-,A)=\prod_{p\in P-\{\hat0\}} \mu(b_{r(p)},a_p) \mvirg$$
where $\mu(b_{r(p)},a_p)$ denotes the M\"obius function for the lattice~$T$.
Now we allow the family~$A$ to vary (i.e. $a_p$ varies in $[b_{r(p)},b_p]$ for each~$p\neq\hat0$) and we define
$$j^B=(-1)^n\sum_{A\in\CF_B}\mu(B^-,A)\, j_A^B \in \Hom_{k\CL}(P,T) \mpoint$$
By Proposition~10.2 of~\cite{BT2}, $f_B=j^B\pi^B$ is an idempotent in~$\End_{k\CL}(T)$
and when $n\geq0$ varies and $B\in \CB_{T,n}$ varies, the idempotents $f_B$ are pairwise orthogonal (Corollary~10.5 of~\cite{BT2}).
This allows us to define the idempotent
$$e_T^{\tot} = \sum_{n=0}^N \sum_{B\in \CB_{T,n}} f_B \mpoint$$
By Theorem 10.8 of~\cite{BT2}, $e_T^{\tot}$ is a central idempotent
and is the identity element of the two-sided ideal $\End_{k\CL}^{\tot}(T)$.
Thus we recover the notation of Section~\ref{Section-idempotent}.\par

\bigskip\noindent
{\bf Second proof of Theorem~\ref{main-theorem}~: }
For each $B\in\CB_{T,n}$, the idempotent $f_B$ is a linear combination of join-morphisms $j_A^B\pi^B$.
We are going to prove that most of these join-morphisms cancel pairwise in the sum
\begin{equation} \label{sum}
e_T^{\tot} = \sum_{n=0}^N \sum_{B\in \CB_{T,n}} (-1)^n\sum_{A\in\CF_B}\mu(B^-,A) \, j_A^B\pi^B \mpoint
\end{equation}
More precisely, we consider all triples $\{(n,B,A)\mid n\in\N, \, B\in\CB_{T,n}, \,A\in\CF_B\}$ such that $a_x<b_x$ for some $x\in P$.
For such a triple, we let $p\in P$ be minimal with respect to the condition $a_p<b_p$.
Since $a_p\in[b_{r(p)},b_p]$ for $p\neq\hat0$, we can either have $b_{r(p)}=a_p< b_p$ or
$b_{r(p)}<a_p<b_p$. The case $p=\hat0$ is special because we always have $a_{\hat0}=\hat0$.
It follows that $p$ must satisfy one of the following 4 cases~:
\par

\begin{enumerate}
\item[{\bf A1.}] $p\neq\hat0$, $r(p)\neq\hat0$, $b_{r(p)}=a_p< b_p$, and $a_x=b_x$ for any $x<p$.
\item[{\bf A2.}] $p\neq\hat0$, $r(p)=\hat0$, $b_{r(p)}=a_p< b_p$, and $\hat0=a_{\hat0}=b_{\hat0}$.
\item[{\bf B1.}] $p\neq\hat0$, $b_{r(p)}<a_p<b_p$, and $a_x=b_x$ for any $x<p$.
\item[{\bf B2.}] $p=\hat0$ and $\hat0=a_{\hat0}<b_{\hat0}$.
\end{enumerate}

{\bf Case A1.} Suppose we are in Case A1. Define
$$\widetilde P=P-\{p\}\,, \qquad\widetilde b_q = b_q \;\;\forall \,q\in \widetilde P-\{r(p)\} \,,\qquad \widetilde b_{r(p)} = b_p \mpoint$$
This defines a chain $\widetilde B$ in~$T$ and a surjective join-morphism $\pi^{\widetilde B}:T\to \widetilde P$, satisfying in particular $\pi^{\widetilde B}(b_p)=r(p)$.
Let $\widetilde A \in \CF_{\widetilde B}$ be the family defined by
$$\widetilde a_q= a_q  \;\;\forall \,q\in \widetilde P-\{r(p)\} \,,\qquad \widetilde a_{r(p)}=b_{r(p)} \mvirg$$
and let $j_{\widetilde A}^{\widetilde B}:\widetilde P\to T$ be the corresponding join-morphism.
Then we obtain
$$b_{r(r(p))} < b_{r(p)} < b_p \,, \qquad\text{ that is, } \qquad \widetilde b_{r(r(p))} <\widetilde a_{r(p)}<\widetilde b_{r(p)} \mvirg$$
so that $\widetilde P$ and its element $r(p)$ are in Case~B1 (because $r(p)\neq\hat0$ by assumption~A1).
Moreover, $j_A^B\pi^B= j_{\widetilde A}^{\widetilde B}\pi^{\widetilde B}$.
This is easy to check on most elements of~$T$, the only nontrivial case being
$$j_A^B\pi^B(b_p) = j_A^B(p)=a_p=b_{r(p)} = \widetilde a_{r(p)}
= j_{\widetilde A}^{\widetilde B}(r(p)) = j_{\widetilde A}^{\widetilde B}\pi^{\widetilde B}(b_p) \mpoint$$
Finally, since $\mu(b_{r(p)},a_p)=\mu(b_{r(p)},b_{r(p)})=1$, the coefficient of $j_A^B\pi^B$ is equal to
\begin{eqnarray*}
(-1)^n \mu(B^-,A)&=&(-1)^n\prod_{x\in P-\{\hat0\}} \mu(b_{r(x)},a_x) \\
&=&(-1)^n\prod_{x\in P-\{\hat0,p\}} \mu(b_{r(x)},a_x) \\
&=&(-1)^n\prod_{x\in \widetilde P-\{\hat0\}} \mu(\widetilde b_{r(x)},\widetilde a_x) \\
&=& -(-1)^{n-1} \mu(\widetilde B^-,\widetilde A) \mvirg
\end{eqnarray*}
using the fact that, for $x=r(p)$, we have $\widetilde a_{r(p)}=b_{r(p)}$ and also $a_{r(p)}=b_{r(p)}$ by minimality of the choice of~$p$.
This shows that
$$(-1)^n \mu(B^-,A) \, j_A^B\pi^B \qquad \text{ and } \qquad (-1)^{n-1} \mu(\widetilde B^-,\widetilde A) \, j_{\widetilde A}^{\widetilde B}\pi^{\widetilde B}$$
cancel in the sum~(\ref{sum}).
Thus any Case~A1 cancels with some Case~B1.\par

\medbreak
{\bf Case B1.} Suppose we are in Case B1. Define
$$\widehat P=P_{< p} \sqcup \{s\} \sqcup P_{\geq p} \mvirg$$
with the total order defined by $x<s$ for all $x\in P_{< p}$ and $s<x$ for all $x\in P_{\geq p}$, so that $r(p)=s$.
Moreover, define
$$\widehat b_q = b_q \;\;\forall \,q\in \widehat P-\{s\} \,,\qquad \widehat b_s = a_p \mpoint$$
This defines a chain $\widehat B$ in~$T$ and a surjective join-morphism
$\pi^{\widehat B}:T\to \widehat P$,
satisfying in particular $\pi^{\widehat B}(a_p)=s$ and $\pi^{\widehat B}(b_p)=p$.
Finally, let $\widehat A \in \CF_{\widehat B}$ be the family defined by
$$\widehat a_q= a_q  \;\;\forall \,q\in \widehat P-\{s\} \,,\qquad \widehat a_s=a_p \mvirg$$
and let $j_{\widehat A}^{\widehat B}:\widehat P\to T$ be the corresponding join-morphism.
Then we obtain
$$a_p < b_p \,, \qquad\text{ that is, } \qquad \widehat b_{r(p)} = \widehat a_p < \widehat b_p \mvirg$$
so that $\widehat P$ and its element $p$ are in Case~A1 (because $r(p)=s\neq\hat0$).
Applying the procedure described in Case~A1, we note that $\widehat P-\{p\}$ is isomorphic to~$P$
and it follows easily that we recover the Case~B1 we started with.
Thus every Case~B1 has been canceled with a corresponding Case~A1.\par

\medbreak
{\bf Case A2.} Suppose we are in Case A2.
Since $r(p)=\hat0$, $p$ is the least element of~$P-\{\hat0\}$.
Define
$$\widetilde P=P-\{p\}\,, \qquad\widetilde b_q = b_q \;\;\forall \,q\in \widetilde P-\{\hat0\} \,,\qquad \widetilde b_{\hat0} = b_p \mpoint$$
This defines a chain $\widetilde B$ in~$T$ and a surjective join-morphism $\pi^{\widetilde B}:T\to \widetilde P$, satisfying in particular $\pi^{\widetilde B}(b_p)=\hat0$.
Let $\widetilde A \in \CF_{\widetilde B}$ be the family defined by
$$\widetilde a_q= a_q  \;\;\forall \,q\in \widetilde P-\{\hat0\} \,,\qquad \widetilde a_{\hat0}=b_{\hat0} \mvirg$$
and let $j_{\widetilde A}^{\widetilde B}:\widetilde P\to T$ be the corresponding join-morphism.
We have $\hat0=b_{\hat0}$ by minimality of~$p$ and we obtain
$$\hat0=b_{\hat0} < b_p \,, \qquad\text{ that is, } \qquad \hat0 = \widetilde a_{\hat0}<\widetilde b_{\hat0} \mvirg$$
so that $\widetilde P$ and its element $\hat0$ are in Case~B2.
The argument for the M\"obius function holds in the same way as in Case~A1 and it follows that any Case~A2 cancels with some Case~B2 in the sum~(\ref{sum}).\par

\medbreak
{\bf Case B2.} Suppose we are in Case B2.
Define
$$\widehat P= \{s\} \sqcup P \mvirg$$
with the total order defined by $s<x$ for all $x\in P$, so that $r(p)=s=\hat0_{\widehat P}$.
Moreover, define
$$\widehat b_q = b_q \;\;\forall \,q\in P \,,\qquad \widehat b_s = \hat0 \mpoint$$
This defines a chain $\widehat B$ in~$T$ and a surjective join-morphism
$\pi^{\widehat B}:T\to \widehat P$,
satisfying in particular $\pi^{\widehat B}(b_p)=p$.
Finally, let $\widehat A \in \CF_{\widehat B}$ be the family defined by
$$\widehat a_q= a_q  \;\;\forall \,q\in P \,,\qquad \widehat a_p=\hat0 \mvirg$$
and let $j_{\widehat A}^{\widehat B}:\widehat P\to T$ be the corresponding join-morphism.
Then we obtain
$$\hat0 < b_p \,, \qquad\text{ that is, } \qquad \hat0=\widehat a_{\hat0}=\widehat b_{\hat0} = \widehat a_p < \widehat b_p \mvirg$$
so that $\widehat P$ and its element $p$ are in Case~A2.
Applying the procedure described in Case~A2, we note that $\widehat P-\{p\}$ is isomorphic to~$P$
and it follows easily that we recover the Case~B2 we started with.
Thus every Case~B2 has been canceled with a corresponding Case~A2.\par

\medbreak
Applying the cancelations described above, we can now eliminate all the join-morphisms $j_A^B\pi^B$
corresponding to a triple $(n\in \N, B\in\CB_{T,n}, A\in \CF_B)$ satisfying $a_x<b_x$ for some $x\in P$.
We are left with the triples satisfying $a_x=b_x$ for all $x\in P$.
In such a case, we have $b_{\hat0}=\hat0$, that is, $\hat0\in B$, hence $B\in\CZ_{T,n}$.
Moreover, $A=B-\{\hat0\}=:B^+$ and
$j_{B^+}^B\pi^B(t)=b_p$ if $\pi^B(t)=p$, that is, if $t\leq b_p$ but $t\not\leq b_{r(p)}$.
In other words,
$$j_{B^+}^B\pi^B(t)=\min\{b\in B\mid t\leq b\}$$
and this is exactly the definition of the endomorphism $\alpha_B$ considered in Section~\ref{Section-idempotent}. Thus
$$j_{B^+}^B\pi^B=\alpha_B \qquad \forall \; B \in\CZ_{T,n} \mpoint$$
Moreover, the coefficient of $j_{B^+}^B\pi^B$ in the expression for~$e_T^{\tot}$ is the M\"obius function
$$(-1)^n \mu(B^-,B^+) = (-1)^n \prod_{p\in P-\{\hat0\}} \mu(b_{r(p)},b_p)=-\mu(B,\infty) \mvirg$$
by Lemma~\ref{Mobius}, where the latter M\"obius function is the M\"obius function of the poset $\CZ_T\sqcup\{\infty\}$.
It follows that the expression for~$e_T^{\tot}$ given in~(\ref{sum}) reduces to
$$e_T^{\tot} = \sum_{n=0}^N \sum_{B\in \CZ_{T,n}} (-1)\mu(B,\infty) \, \alpha_B \mpoint$$
This completes the second proof of Theorem~\ref{main-theorem}.

\result{Remark} \label{semi-simple}
{\rm 
It is proved in Theorem~10.6 of~\cite{BT2} that the two-sided ideal $\End_{k\CL}^{\tot}(T)$ is isomorphic to a direct sum of matrix algebras
$$\End_{k\CL}^{\tot}(T) \cong \bigoplus_{n=0}^N M_{|\CZ_{T,n}|} (k) \mvirg$$
where $N$ is the maximal length of a chain in~$T$.
It should be noticed that the new approach to the idempotent $e_T^{\tot}$ explained in the present paper does not simplify in any way the proof of this result.
In particular, if $T$ is totally ordered, then
$$\End_{k\CL}(T)=\End_{k\CL}^{\tot}(T)\cong \bigoplus_{n=0}^N M_{|\CZ_{T,n}|} (k) \mvirg$$
and this is a semi-simple algebra whenever $k$ is a field.
As noticed in Remark~11.3 of~\cite{BT2}, this result is similar, but not equivalent, to a theorem proved in~\cite{FHH} about the planar rook algebra.
}
\fresult


\bigskip
\noindent
Serge Bouc, CNRS-LAMFA, Universit\'e de Picardie - Jules Verne,\\
33, rue St Leu, F-80039 Amiens Cedex~1, France.\\
{\tt serge.bouc@u-picardie.fr}

\medskip
\noindent
Jacques Th\'evenaz, Section de math\'ematiques, EPFL, \\
Station~8, CH-1015 Lausanne, Switzerland.\\
{\tt Jacques.Thevenaz@epfl.ch}


\begin{thebibliography}{}

\bibitem[BT1]{BT1}
S.~Bouc, J.~Th\'evenaz.
\newblock Correspondence functors and finiteness conditions,
\newblock {\em J. Algebra} 495 (2018), 150--198.

\bibitem[BT2]{BT2}
S.~Bouc, J.~Th\'evenaz.
\newblock Correspondence functors and lattices,
\newblock preprint, 2017.

\bibitem[BT3]{BT3}
S.~Bouc, J.~Th\'evenaz.
\newblock The algebra of Boolean matrices, correspondence functors, and simplicity,
\newblock in preparation.

\bibitem[FHH]{FHH}
D.~Flath, T.~Halverson, K.~Herbig.
\newblock The planar rook algebra and Pascal's triangle,
\newblock {\em Enseign. Math.} 55 (2009), no.~1-2, 77--92.

\bibitem[Qu]{Qu}
D.~Quillen,
\newblock Homotopy properties of the poset of nontrivial $p$-subgroups of a group,
\newblock {\em Adv. Math.} {\bf 28} (1978), 101--128.

\bibitem[St]{St}
R.~P.~Stanley.
\newblock {\em Enumerative Combinatorics,  Vol.~I}, Second edition,
\newblock Cambridge studies in advanced mathematics 49,
\newblock Cambridge University Press, 2012.


\end{thebibliography}
\end{document}